\documentstyle{amsppt}
\magnification = 1200
\hcorrection{.25in}

\document
\topmatter 

\title  {Divergence of projective structures and lengths of measured
laminations} 
\endtitle 

\rightheadtext {Divergence of projective structures}

\author Harumi Tanigawa  
\endauthor  

\address Graduate School of Polymathematics,  Nagoya University,
Nagoya 464-01,  Japan 
\endaddress

\abstract  Given a complex structure, we investigate diverging sequences
of projective structures on the fixed complex structure in terms of
Thurston's parametrization.  In particular, we will give a geometric proof
to the theorem by Kapovich stating that as the
projective structures on a fixed complex structure diverge so do their
monodromies. In course of arguments, we extend the concept of realization
of laminations for
$\operatorname {PSL} (2, \bold C)$-representations of surface groups.

\endabstract

\endsubjclass
\thanks {Research at MSRI is supported by NSF grant 
\#DMS--9022140}
\endthanks
\endtopmatter

\head 1. Introduction
\endhead

Throughout  this paper, $\Sigma_g$ denotes an oriented differentiable
surface of genus $g >1$ and $P_g$ the space of $\operatorname {\bold
CP}^1$-structures (or, projective structures) on $\Sigma_g$.  Here, a
$\operatorname {\bold CP}^1$-structure on a surface is a structure
modelled on 
$(\operatorname {\bold CP}^1, \operatorname {PSL}(2,{\bold C}))$. 
(Here $\operatorname {PSL}(2,{\bold C}$) is identified with the
projective automorphism group of
$\operatorname {\bold CP}^1$.) As a M\"obius transformation is
holomorphic, each $\operatorname {\bold CP}^1$-structure determines
its underlying complex structure.\par
 There are two parametrizations of
$P_g$ which are often quoted.
 One of them is by the bundle
$\pi: Q_g \to T_g$ of holomorphic quadratic differentials on Riemann
surfaces over Teichm\"uller space $T_g$. The canonical projection
$\pi$ represents the correspondence from each $\operatorname {\bold
CP}^1$-structure to its underlying complex structure. As for the
detail about this parametrization,  see Hejhal
\cite {H}, for example.
 The other parametrization was introduced by Thurston, that says
every $\operatorname {\bold CP}^1$-structure is obtained from some 
hyperbolic structure by grafting a measured lamination to it:
$P_g$ is parametrized  by
$T_g
\times
\Cal {ML}$, where
$\Cal {ML}$ stands for the space of measured laminations. As for a
description of grafting and this parametrization, see
\cite {KT} or \cite {Tg}. 
\par

In this paper, we will investigate 
$\operatorname {\bold CP}^1$-structures on a fixed complex structure in
terms of Thurston's parametrization. In particular, we will give a
geometric proof for the fact, originally proved by Kapovich
\cite {Ka}, that the monodromies diverge as the $\operatorname {\bold
CP}^1$-structures diverge on a fixed complex structure.

\par  The idea  is to generalize the concept of realization of
measured laminations in hyperbolic $3$-manifolds to that for
$\operatorname {PSL}(2, {\bold C})$-representations of the
fundamental group
$\pi_1(\Sigma_g)$ which is not necessarily discrete, and then 
observe that the lengths of measured laminations are locally bounded in
the space of $\operatorname {PSL}(2,\bold C)$-representations, while the
length of some measured lamination diverges as 
$\operatorname {\bold CP}^1$-structures on the fixed complex
structure diverge.

\par The author would like to thank William Thurston for inspiring
conversations. This work was done at Mathematical Sciences Research
Institute. The author is grateful to their hospitality.

\head  2. Pleated surfaces and measured laminations for
representations
\endhead 
 First, we recall the monodromy of a
 $\operatorname {\bold CP}^1$-structure. Given a
$\operatorname {\bold CP}^1$-structure on
$\Sigma_g$, we take its {\it developing map}: beginning with a local
chart of the
$\operatorname {\bold CP}^1$-structure, we take its
 analytic continuation along curves and have a multivalued
meromorphic function (meromorphic with respect to the underlying
complex structure of the 
$\operatorname {\bold CP}^1$-structure). The values over a point
differ from each other by post-compositions of M\"obius
transformations, which are determined by the closed curves along which
the analytic continuation is carried out. Thus this procedure yields 
a homomorphism 
$\pi_1 \Sigma_g \to \operatorname {PSL}(2,{\bold C})$, which is
determined up to conjugations by M\"obius transformations (depending
on the local chart to begin the analytic continuation with). This
homomorphism is called a  {\it monodromy} of the $\operatorname {\bold
CP}^1$-structure, or the {\it monodromy representation of }
$\pi_1 \Sigma_g$.\par Now we review quickly the two parametrization of
$P_g$ mentioned in the introduction. \par First we recall the
parametrization by the bundle of quadratic differentials
$\pi : Q_g \to T_g$. Given a $\operatorname {\bold CP}^1$-structure,
we take its developing map $f$. Then its Schwarzian derivative,
defined by
$(f''/f')^2 - 1/2 (f''/f')'$, is a quadratic differential holomorphic
with respect to the underlying complex structure. (For geometric
meanings of Schwarzian derivatives see Thurston \cite {Th2}.)
Therefore, we have a Riemann surface $X \in T_g$ and a holomorphic
quadratic differential on $X$. Conversely, it is known that given a
Riemann surface $X$ and a holomorphic quadratic differential $q$
there is a $\operatorname {\bold CP}^1$-structure on $X$ whose
Schwarzian derivative of the developing map is $q$. Note that for
each point $X \in T_g$, the space of $\operatorname {\bold
CP}^1$-structures on $X$ is the fiber $Q_g(X)$ of $\pi : Q_g \to T_g$
over
$X$, which is complex $3g-3$-dimensional Banach space.
\par
On the other hand, Thurston's parametrization theorem is as follows.
\proclaim {Theorem 2.1 (Thurston, unpublished)} The space of
$\operatorname {\bold CP}^1$-structures on $\Sigma_g$ is parametrized
by $T_g \times \Cal {ML}$: every $\operatorname {\bold
CP}^1$-structure is obtained by grafting a measured lamination to a
hyperbolic surface.
\endproclaim 
As for a full proof of this parametrization theorem, see section 2 of \cite
{KT},or, for a short discription,  section 2 of \cite {Tg}.
His idea is as follows: Given a hyperbolic surface $X$ and a measured
geodesic  lamination
$\lambda$, first, embed the universal cover $\operatorname {\bold H}^2$ of $X$
into $\operatorname {\bold H}^3$. Then, lift $\lambda$ to $\operatorname {\bold H}^2$
and then bend $\operatorname {\bold H}^2$ along the lift of $\lambda$ so that the
bending measure is the measure of $\lambda$. This bent structure can be ``pushed
forward" to the sphere at infinity to determin a $\operatorname {\bold
CP}^1$-structure.
\par
Note that the projection from $T_g \times \Cal {ML}$ to its
first coordinate $T_g$ is {\it not} the correspondence between
$\operatorname {\bold CP}^1$-structures and their underlying complex
structures. On the contrary, the underlying complex structure of any
$\operatorname {\bold CP}^1$-structure represented by $(X , \mu) \in T_g
\times \Cal {ML}$, differs from $X$, unless
$\mu = 0$.
\remark {Remark} In relation to this parametrization,
Labourie \cite {L} gave a new parametrization for $P_g$.
\endremark 
\par

Now we define the pleated maps and the lengths of measured
laminations for 
$\operatorname {PSL}(2,{\bold C})$-representations of $\pi_1
\Sigma_g$.
  We begin with ``pleated surfaces". It is defined in a parallel way
to that of
$3$-manifolds.  .

\definition {Definition 2.2} Let $\phi : \pi_1 \Sigma_g \to
\operatorname {PSL}(2,{\bold C})$ be a representation.
 A {\it pleated mapping} for $\phi$ is a continuous mapping 
$f : \operatorname {\bold H}^2 \to \operatorname {\bold H}^3$ with
the following properties:

\roster
\item  there exists a Fuchsian group $\Gamma$ acting on
$\operatorname {\bold H}^2$ isomorphic
 to $\pi_1 \Sigma_g$ such that $f \circ \gamma = \phi (\gamma) \circ
f $ for all $\gamma \in \Gamma$.
\item for every $z \in \operatorname {\bold H}^2$, there is an open
interval of hyperbolic line containing $z$ which is mapped by $f$ to
a straight line segment.
\item $f$ maps every geodesic segment in $\operatorname {\bold H}^2$ to a
rectifiable arc in $\operatorname {\bold H}^3$ with the same length.
\endroster
\enddefinition 
This is the same as the usual definition of pleated
surfaces when  $\phi (\pi_1 \Sigma_g)$ is discrete. We shall call the
image of a pleated map a {\it pleated surface}.\par

A typical example  of a pleated surface is the ``bent surface" defining
a $\operatorname {\bold CP}^1$-structure in Thurston's geometric parametrization
theorem. It is a pleated map for
the monodromy of the $\operatorname {\bold CP}^1$-structure.(see
\cite {EM}, \cite {KT} and \cite {Tg}). \par

 Now we define the lengths of measured laminations for $\operatorname
{PSL}(2, {\bold C})$-representations of $\pi_1 \Sigma_g$.

\definition  {Definition 2.3} Let $\phi : \pi_1 \Sigma_g \to 
\operatorname {PSL}(2, {\bold C})$ be a non-elementary representation
and let $f : \operatorname {\bold H}^2 \to
\operatorname {\bold H}^3$ be  a pleated map  equivariant with respect  to
some Fuchsian group $\Gamma$ isotopic to $\pi_1
\Sigma_g$ acting on $\operatorname {\bold H}^2$ and $\phi (\pi_1
\Sigma_g)$. Such a mapping determines a hyperbolic structure on
$\Sigma_g$. Then let
$l_f(\mu)$ be  the total mass of the product of the transverse
measure of
$\mu$ and the length along the lines of the support of $\mu$. The
length
$l_\phi(\mu)$ is the infimum of $l_f(\mu)$ where the infimum is taken
over all  Fuchsian groups $\Gamma$ isotopic to $\pi_1
\Sigma_g$ and pleated maps $f$.
\enddefinition 
\remark {Remark} As we will see in Lemma 2.5, for any non-elementary 
$\operatorname {PSL (2, \bold C)}$-rep\-res\-en\-ta\-tion of $\Sigma_g$, there
exist at least one pleated map as above.
\endremark  As Definition 2.2, when the
$\operatorname {PSL}(2, {\bold C})$-representation  of $\pi_1 \Sigma_g$ is
discrete, the above definition coincides with that of usual length of a
measured lamination. Also, given a homotopy class of an equivariant map,
the uniqueness of the realization of a lamination holds by the same reason
as in the discrete case.
\par
 The length function is a function of two variables:
 measured laminations and representations.  We might as well expect the
continuity of the length function with respect to the both variables.  (In
fact, for a simple closed curve, the length is continuous with respect to 
$\operatorname {PSL}(2, {\bold C})$-representation of
$\pi_1 \Sigma_g $.)  However, we do not need that strong property for
our purpose. What we need is local boundedness (Lemma 2.5 bellow). \par 
 Kapovich (\cite {Ka}) has already mentioned the concept of pleated
maps and showed the existence for any non-elementary 
$\operatorname {SL}(2,{\bold C})$-representation of $\pi_1 \Sigma_g$
briefly. We use these
  examples of pleated maps observed by Kapovich.\par Here, for
convenience, we will exibit somewhat detailed argument about them.
\par

The following decomposition is the key.
\proclaim {Theorem 2.4 (\cite {Ka})} For any non-elementary
representation 
$$\phi : \pi_1\Sigma_g  \to
\operatorname {PSL}(2,{\bold C})$$ there exists a pants decomposition 
$\Sigma_g = P_1 \cup ...\cup P_{2g-2}$ such that 
\roster 
\item $\phi (\pi_1(P_i))$ is non-elementary and
\item $\phi (\gamma)$ is loxodromic for each boundary curve $\gamma$
of $P_i$  for $ i = 1,...,2g-2$.
\endroster
\endproclaim 

\proclaim {Lemma 2.5} The length function of measured laminations is
locally bounded in the space of non-elementary $\operatorname {PSL}(2,\bold
C)$-representations of $\pi_1\Sigma_g$.
\endproclaim
\demo {Proof}  Given a non-elementary representation $\phi :
\pi_1\Sigma_g  \to
\operatorname {PSL}(2,{\bold C})$ take a pants decomposition 
$\Sigma_g = P_1 \cup ...\cup P_{2g-2}$ as in the above theorem.
Decompose each pants into two ideal triangles by adding leaves
spiraling around the boundaries. \par Take a continuous equivariant
map $f : \operatorname {\bold H}^2 \to
\operatorname {\bold H}^3$ (with respect to
 any Fuchsian group isotopic to $\Sigma_g$ and $\phi (\Sigma_g)$) 
and homotope $f$ so that $f$ sends each leaf of the ideal
triangulation onto a straight line, then straighten the map in each
ideal triangle. This is possible because by Theorem 2.4 the
restriction of $\phi$ to the fundamental group of each pants is
isomorphic. Thus we have a pleated surface with pleating locus the
geodesics defining the ideal triangulation.\par Thus, for each
non-elementary
$\operatorname {PSL}(2, {\bold C})$-representation $\phi$, 
the length of
any measured lamination is finite. Here, note that when we take an ideal
triangulation as above for $\phi$, we can take the same one
 for
$\operatorname {PSL}(2, {\bold C})$-representations which are sufficiently
close to
$\phi$. As the $\operatorname {PSL}(2, {\bold C})$-representations
changing continuously near
$\phi$, we can take  hyperbolic structures on $\Sigma_g$
determined by the pleated maps with the fixed ideal triangulation so that they
change continuously. From the continuity of the lengths of measured
laminations on surfaces, we have a local upper bound of the lengths of measured
laminations by the lengths on the hyperbolic structures. It follows that the lengths
for the representations are locally bounded.
\qed
 \enddemo
\remark {Remark} Let $mon : P_g
\to \operatorname {Hom} (\Sigma_g, 
\operatorname {PSL}(2, {\bold C}))$ be the monodromy map, namely, the
map
$mon$ sends each $\operatorname {\bold CP}^1$-structure to its
monodromy representation. Hejhal \cite {H} showed that $mon$ is a
local diffeomorphism. Therefore, the image $mon(P_g)$ is a region of
$\operatorname {Hom} (\Sigma_g, 
\operatorname {PSL}(2, {\bold C}))$. Therefore, obviously the lengths
of measured laminations are locally bounded on the region $mon (P_g)$ by
Theorem 2.1.  A theorem by Gallo, Kapovich and Marden \cite {GKM} says
that  every non-elementary
$\operatorname {SL}(2, {\bold C})$-representation of
$\Sigma_g$ is obtained as the monodromy of a
$\operatorname {\bold CP}^1$-structure.  Therefore, if we would use that
result, the local boundedness on the entire space of $\operatorname
{\bold CP^1}$-structures would follow immediately. However, 
we exibited a direct proof rather than using their theorem.
\endremark

\head
3. Divergence of projective structures
\endhead

Our goal is to give a geometric proof for the following:
\proclaim {Theorem 3.1 (\cite {Ka})} For any compact set $K$ of $T_g$
and any diverging sequence $\{q_n\}$ of $\operatorname {\bold
CP}^1$-structures with $\{\pi (q_n)\} \subset K$, their monodromies
$mon (q_n)$ diverge. Here, $\pi(q_n)$ stands for the underlying
complex structure of $q_n$. 
\endproclaim

The idea is as follows: The extremal length  is the square order of the hyperbolic
length when hyperbolic structures stay in a compact set of $T_g$.
On the other hand,
Theorem 3.2 bellow will show that
$l_X(\mu)$ and
$E_{gr_\mu(X)}(\mu)$ is in the same order when they are large.
Here, $gr_\mu(X)$ denote the Riemann surface obtained by grafting $\mu$ to
$X$ (namely, the underlying structure of the $\operatorname {\bold
CP}^1$-structure given by  Thurston's parameter $(X, \mu)$), 
$l_X(\mu)$ is the hyperbolic length of
$\mu$ on $X$ and $E_{gr_\mu (X)}(\mu)$ is the extremal length of
$\mu$ on the grafted surface
$gr_\mu(X)$. Furthermore, by the definition of bendinging, $l_X(\mu)$
is the length of $\mu$ for the monodromy of the $\operatorname {\bold
CP}^1$-structure.
 Therefore, for a diverging sequence $q_n = (X_n, \mu_n)$  as in the assumption of
Theorem 3.1, 
$l_{X_n}(\mu)$ had to grow ``too  fast" and we can not maintain
the monodromy in a bounded set of
$\operatorname {PSL}(2, \bold C)$-representation as $n \to \infty$,
considering Lemma 2.5.\par

\proclaim {Theorem 3.2(\cite {Tg})}
 Let $X$ be a hyperbolic surface and
$\mu$ be a measured lamination. Let
$h : gr_\mu(X) \to X$ be the harmonic map with respect to the
hyperbolic metric on $X$ and $\Cal {E} (h)$ be its energy.

Then
$$\frac {1}{2} l_X(\mu) \le \frac {1}{2}\frac{l_X(\mu)^2}{E_{gr_\mu
(X)}(\mu)} \le 
\Cal {E} (h) \le \frac {1}{2} l_X(\mu) + 4\pi (g-1).$$
\endproclaim 
Recall that in Remark 1 in \cite {Tg}, we observed that
when
$l_X(\mu)$ is very large, that is, when the grafted structure is far
from $X$, grafting is very close to the inverse of the harmonic map
$h$ and the pleated surface defining the grafting is close to the
image of the harmonic equivariant map for the monodromy
representation. \par

The following fact follows easily from the above inequality (cf. \cite
{Tg}).

\proclaim {Corollary 3.3} For any $X \in T_g$, the mapping $gr_\cdot
(X) :
\Cal {ML} \to T_g$ is a proper map. For any $\mu \in \Cal {ML}$, the
mapping $gr_\mu(\cdot) : T_g \to T_g $ is a proper map.
\endproclaim

\demo {Proof}
The properness of $gr_\mu(\cdot): T_g \to T_g $ was shown in \cite
{Tg}. For the properness of $gr_\cdot (X) : \Cal {ML} \to T_g$, note
that

$$ 
  \frac {E_{gr_\mu (X)}(\mu)}{E_X (\mu)} \le
          \frac {l_X(\mu)}{E_X (\mu)}, $$
by Theorem 3.2.
When $X$ is fixed and $\mu$ tends to infinity, the right term tends
to $0$, therefore, so is the left term. It follows that 
the Riemann surface $gr_\mu (X)$ tends to infinity.\qed

\enddemo
\remark {Remark} In fact, in what follows we   use only the leftmost
inequality in Theorem 3.2. The others was necessary to show that 
$gr_\mu(\cdot)$ is proper.
\endremark

Now we  can present a geometric proof of Theorem 3.1.
\demo{Proof of Theorem 3.1}  

 Let $(X_n, \mu_n) \in T_g \times \Cal {ML}$ be the Thurston
coordinate of the projective structure $q_n$ and denote the underlying complex
structure by 
$\pi (q_n) = gr_{\mu_n}(X_n)$. Note that in Corollary 3.3  in fact we can show that
the mappings are ``locally uniformly" proper with respect to each variable.
Therefore, when
$\pi(q_n)$ stays in $K$ and $\{q_n\}$ tends to infinity,
 both of
$\{X_n\}$ and $\{\mu_n\}$ tend to infinity. Take a sequence $\epsilon_n$ converging to $0$ such that $\epsilon_n \mu_n$
converges to a non-zero measured lamination $\mu$, taking a subsequence
of $\{\mu_n\}$ if necessary.
\par
 By
Theorem 3.1, 
 $$ \frac{l_{X_n}(\mu_n)^2} {l_{X_n}(\mu_n) + 8\pi (g-1)} 
      \le E_{gr_{\mu_n}(X_n)}(\mu_n)  \le  l_{X_n}(\mu_n). $$  
Multiply
each term of the above inequality by
$\epsilon_n^2$. Then we have 
   $$ \frac{l_{X_n}({\epsilon_n \mu_n})^2} {l_{X_n}({\mu_n}) + 8\pi
(g-1)} 
      \le E_{gr_{\mu_n}(X_n)}(\epsilon_n\mu_n)  \le 
\epsilon_n l_{X_n}(\epsilon_n \mu_n). $$

Note that $l_{X_n}(\epsilon_n \mu_n)$ is the length of 
$\epsilon_n \mu_n$ for the monodromy representation $mon(q_n)$ because the
pleated surface defining the Thurston coordinate for $q_n =
(X_n, \mu_n)$ has pleating locus
$\mu_n$. Therefore, if the representation $mon(q_n)$ did not diverge, 
$l_{X_n}(\epsilon_n \mu_n)$ would converge  to a non-negative number by
Lemma 2.5, taking a subsequence if necessary. Therefore, by the above
inequality,
$$\lim_{n \to \infty} E_{gr_{\mu_n}(X_n)}(\epsilon_n\mu_n) = 0$$
 However, by taking a subsequence if necessary, we may assume that
$\pi(q_n) = gr_{\mu_n}(X_n)$ converges to a point $Y \in T_g$ by the
assumption that $\{q_n\}$ stays in the compact set $K$. Therefore
$E_{gr_{\mu_n}(X_n)}(\epsilon_n\mu_n)$ converges to 
$E_Y(\mu)$, which is a positive number. This is a contradiction.
\qed

\enddemo

\enddocument